\documentclass[10pt,a4paper]{article}
\usepackage[latin1]{inputenc}
\usepackage[T1]{fontenc}
\usepackage{amsmath,amsfonts,amssymb}
\usepackage{mathrsfs}
\usepackage[english]{babel}
\usepackage{enumitem}

\usepackage{multicol}
\usepackage[left=3cm,right=3cm,top=2cm,bottom=2cm]{geometry}

\newcommand{\N}{\mathbb{N}}
\newcommand{\Z}{\mathbb{Z}}

\newcommand{\R}{\mathbb{R}}

\let\iff=\Longleftrightarrow

\newcommand{\Erdos}{Erd\H{o}s-Straus}
\newcommand{\et}{\mathrm{\quad and \quad}}
\newcommand{\ou}{\mathrm{\quad or \quad}}

\newcommand{\ppcm}{\mathop{\rm LCM}}
\newcommand{\Hline}{\hline \vphantom{\Big(}}
\newcommand{\dg}{\ensuremath{d^{\,\circ}}} 

\newenvironment{enum}{\begin{enumerate}[label=$\bullet$, topsep=0.75em, itemsep=0.5em, labelsep=1ex]}{\end{enumerate}}

\newenvironment{dem}{\textsc{Proof}  }{\hfill $\blacksquare$ \vskip 5mm plus 1mm minus 2mm}

\newtheorem{lemme}{\sc Lemma}
\newtheorem{prop}{\sc Proposition}
\newtheorem{corolaire}{\sc Corollary}

\begin{document}
\author{Serge E. \textsc{Salez}}
\title{The \Erdos\ conjecture\\
New modular equations \\
and checking up to $N=10^{17}$}
\date{}
\maketitle

\parindent=0pt

\begin{abstract}
In 1999 Allan Swett \cite{Swett} checked (in 150 hours) the  \Erdos\ conjecture up to $N=10^{14}$ with a sieve based on a single modular equation. After having proved the existence of a  "complete" set of seven modular equations (including three new ones), this paper offers an optimized sieve based on these equations. A program written in  C++ (and given elsewhere) allows then to make a checking whose running time, on a typical computer%
\footnote{AMD TurionII Dual-Core Mobile M250 ($64$ bits, $16\,100$ MIPS).}, %
range from few minutes for $N=10^{14}$ to about 16 hours for $N=10^{17}$. 
\end{abstract}

\section{Basic formulas}

A fraction is said to be \textit{$k$-Egyptian} if it is the sum of at most $k$ positive unit fractions (i.e with numerator equal to 1). The \Erdos\ conjecture claims that $4/n$ is a 3-Egyptian fraction for any $n>1$.

\subsection{Reduction} \label{reduction} 

Through the identities
\[ \frac{1}{t}=\frac{1}{t+1}+\frac{1}{t(t+1)} \]

\[ \frac{2}{2t-1}=\frac{1}{t}+\frac{1}{t(2t-1)} \]
it is  equivalent (for $n>2$) to require having \textit{exactly} 3 \textit{different} unit fractions,  what we shall do thereafter.\\ 

On the other hand, the identities    

\[ \frac{4}{3t-1}=\frac{1}{t}+\frac{1}{3t-1}+\frac{1}{t(3t-1)} \]

\[ \frac{4}{4t-1}=\frac{1}{t}+ \frac{1}{t(4t-1)} \]

\[ \frac{4}{8t-3}=\frac{1}{2t}+\frac{1}{t(8t-3)}+\frac{1}{2t(8t-3)}\]

show that the conjecture is verified  if $n=-1\mod 3$ or $n=-1 \mod 4$ or $n=-3 \mod 8$. Moreover, if $4/n$ is 3-Egyptian then $4/kn$ is too.\\

To conclude, it is then sufficient to prove that $4/p$ is 3-Egyptian, for any prime integer $p$ such that $p=1 \mod 24$.

\subsection{Rosati's formulas}

The following proposition is due (according to Mordell\footnote{%
So, unlike most paper, we don't attribute to Mordell what Mordell himself attribute to others mathematicians. In his book \cite{Mordell}, often quoted, the four pages given to this conjecture doesn't introduce a personal work but report briefly some papers whose sources are scrupulously pointed out : hence, it is absolutely incorrect to speak of "Mordell's theorem" or of "Mordell's formulas". On a different scale, it should be better not to remake what was done with Pell's equation.})
to Rosati \cite{Rosati}. The proof needs only simple calculations and has been given many a time. This one is nevertheless original and standardize the notations.\\

We set $\mathcal{A}=\Z$. Let $\mathcal{A}_+$ be the set of strictly positive elements of $\Z$. In this context, we call \textit{prime element} an odd prime integer.

\begin{prop}  \label{Rosati}
Let $p$ a prime element. The fraction $4/p$ is 3-Egyptian if and only if there exists four elements of $\mathcal{A}_+$ denoted by $A$, $B$, $C$, $D$ such that
\begin{equation}  \label{Rosati_1} 
                 4ABCD=A+B+pC  \et (ABD,p)=1
\end{equation}
or 
\begin{equation}  \label{Rosati_2}
                 4ABCD=p(A+B)+C  \et (ABCD,p)=1
\end{equation}
\end{prop}

\begin{dem}
If we assume that $4/p$ is 3-Egyptian, then there exists 3 elements of $\mathcal{A}_+$ denoted by $X_1,X_2,X_3$ such that 
\begin{equation}\label{Z} 
   \dfrac{4}{p}=\dfrac{1}{X_1}+\dfrac{1}{X_2}+\dfrac{1}{X_3} 
\end{equation}

The $X_i$ are not all divisible by $p$ for otherwise we would have
\[ 4=\dfrac{1}{X_1/p}+\dfrac{1}{X_2/p}+\dfrac{1}{X_3/p}\leqslant 3 \]

In view of (\ref{Z}) it follows that  
\[ 4X_1X_2X_3=p(X_2X_3+X_3X_1+X_1X_2) \]
which shows that $p$ divides at least one of the $X_i$.
Hence we may set $x_i=X_i/p_i$ where
\[ p_1=p_2=1,\quad p_3=p \et (x_1x_2,p)=1\]
or
\[ p_1=p_2=p,\quad p_3=1 \et (x_3,p)=1\]
depending on $p$ divides exactly one or two $X_i$.\\

Therefore, since $p_2p_3=p$
  \[ 4p_1p\;x_1x_2x_3=p(p_2p_3x_2x_3+p_3p_1x_3x_1+p_1p_2x_1x_2) \]
and hence
  \[ 4x_1x_2x_3=p_3x_2x_3+p_3x_3x_1+p_2x_1x_2 \]
We set
 $D=(x_1,x_2,x_3)$ and $x'_i=x_i/D$. Then
  \[ 4Dx'_1x'_2x'_3=p_3x'_2x'_3+p_3x'_3x'_1+p_2x'_1x'_2 \]
  
At last we set 
     \[ A=(x'_2,x'_3),\quad B=(x'_3,x'_1),\quad C=(x'_1,x'_2) \]
     
Since $(x'_1,x'_2,x'_3)=1$ it follows that $A,B,C$ are pairwise relatively prime. So, we may write 
         \[x'_1=BCt_1,\quad x'_2=CAt_2,\quad x'_3=ABt_3\]
where $t_i\in\mathcal{A}_+$ are pairwise relatively prime. We note that 
		\[(t_1,A)=(t_2,B)=(t_3,C)=1  \]
With these notations, we have
\[ 4D\,BC\,CA\,AB\,t_1t_2t_3=p_3CA\,AB\,t_2t_3 +  p_3AB\,BC\,t_3t_1 + p_2BC\,CA\,t_1t_2  \]
and hence  
\[ 4ABCDt_1t_2t_3=p_3 At_2t_3 +  p_3 Bt_3t_1 + p_2 Ct_1t_2  \]

It follows that $t_1\mid p_3At_2t_3$  which reduce to $t_1\mid p_3$ and hence $t_1=1$ for $(x_1,p_3)=1$. Similar arguments lead to $t_2=1$ and $t_3=1$. Finally
\begin{equation} \label{Rosati_eq}
  4ABCD=p_3A+p_3B+p_2C
\end{equation}

Conversely, if we assume that $A,B,C,D$ verify (\ref{Rosati_eq}), we divide by $pABCD$ and then
 
\[ \frac{4}{p}=\frac{1}{p_1BCD}+\frac{1}{p_2ACD}+\frac{1}{p_3ABD} \]
which shows that $4/p$ is 3-Egyptian%
\footnote{For all purpose, we may write $x=T/A$, $y=T/B$, $z=T/C$ where $T=ABCD$.}.
\end{dem}

We observe that if $p$ is not prime, (\ref{Rosati_eq}) is still \textit{sufficient} but no more \textit{necessary.}

\subsection{Notations}
Henceforth, we systematically make use of the notations of  Proposition~\ref{Rosati}.
We add also $E\in\mathcal A_+$ and $F\in\mathcal A_+$ as follow.\\

By (\ref{Rosati_eq}) and since $(C,p_3)=1$, we have $C\mid A+B$. If we write $E=(A+B)/C$ then $E\in\mathcal A_+$ and (\ref{Rosati_eq}) is equivalent to 

\begin{equation} \label{Rosati_E}
 \begin{cases} 
	4ABD=p_3 E + p_2\\ 
	A+B=CE
\end{cases}
\end{equation}
 
The relation (\ref{Rosati_eq}) may be rewritten $(4BCD-p_3)A=p_3 B+p_2 C $. We set $F=4BCD-p_3$  and then (\ref{Rosati_eq}) is equivalent to
\begin{equation} \label{Rosati_F}
 \begin{cases} 
	FA=p_3 B+p_2 C2\\ 
	F+p_3=4BCD
\end{cases}
\end{equation}
The second equation of (\ref{Rosati_F}) shows that $F\in\mathcal A$, the first one that $F\in\mathcal A_+$.\\
 
Moreover, by (\ref{Rosati_E}) we have $ 4(CE-B)BD=p_3 E + p_2$ and then $ (4BCD-p_3)E=4B^2D+p_2 $. Whence
\begin{equation} \label{Rosati_FE}
 	FE=4B^2D+p_2
\end{equation}

\section{Generalization}

\subsection{Definitions} 
Like for the integers, we say that a rational fraction is $k$-Egyptian if it is the sum of at most $k$ inverses of polynomials of $\Z[X]$.\\

We set $\mathcal{A}=\Z[X]$. Let $\mathcal{A}_+$ be the set of polynomials of $\Z[X]$ whose leading coefficient is strictly positive. In this context, we call \textit{prime element} an  irreducible polynomial of $\mathcal{A}_+$.\\
 
In the ring $\mathcal{A}$, the fundamental theorem of arithmetic is true and the GCD is unique if we request it has to be in $\mathcal{A}_+$. Hence, the Proposition~\ref{Rosati} holds also in this new context, without any change neither in the text nor in the proof. It is the same for $E$ and $F$ as well as the related equations.
  
\subsection{First application}

\begin{prop} \textnormal{(Schinzel's Theorem \label{Schinzel})}\\
Let $a>0$ and $b$ such as $(a,b)=1$. If $4/(at+b)$ is 3-Egyptian%
\footnote{$at+b$ is supposed to be a polynomial (abuse of notation).}
then $b$ is a quadratic non residue modulo $a$.
\end{prop}

\begin{dem}
We write $p(t)=at+b$. There exists $\tau$ such as the polynomials $p,A,B,C,D,E$ take strictly positive values for $t>\tau$.\\

Depending on the case, the equation  $ 4(B-CE)BD=p_3 E + p_2$ may be written   
  \[ 4(CE-B)BD=E+p  \ou 4(CE-B)BD=pE+1\]
So we have
  \[ p=(4BCD-1)E-4B^2D  \ou  pE^2=(4BCDE-1)E-4B^2DE\]
and then, if we write $D'=DE$
  \[ p=-4B^2D \mod 4BCD-1 \ou  pE^2=-4B^2D' \mod 4BCD'-1\]
  
If $b$ is a quadratic residue modulo $a$, there exists an integer $k>\tau$ such that $ak+b$ is a square.
If $t=k$, it follows from propriety of the Jacobi symbol%
\footnote{The same notations  $p,A,B,C,D,E$ are still used for the values at $t=k$ of these polynomials. By the way, a similar calculation using the Kronecker symbol is made in the paper of Yamamoto\cite{Yamamoto}.}
  \[ \left( \dfrac{p}{4BCD-1} \right) =  \left( \dfrac{-4B^2D}{4BCD-1} \right) =-1\]
which contradicts the fact that $p$ is a square. Idem with $pE^2$.\\ 

More precisely, if we write $D=2^\alpha m$ where $m$ is odd, we obtain
    
 \[ \left( \dfrac{-4B^2D}{4BCD-1} \right)=-  \left( \dfrac{D}{4BCD-1} \right) =-\left( \dfrac{2}{4BCD-1}\right)^\alpha \left( \dfrac{m}{4BCD-1}\right) \]
 
If $\alpha >0$ then $4ABCD-1=7 \mod 8$ and this implies
 \[\left( \dfrac{2}{4BCD-1}\right)=1 \]   
For the second factor, using the law of quadratic reciprocity, we have  
 \[ \left( \dfrac{m}{4BCD-1} \right)=(-1)^{(m-1)/2}\left( \dfrac{4BCD-1}{m} \right)=(-1)^{(m-1)/2}\left( \dfrac{-1}{m} \right) \]
and then
 \[ \left( \dfrac{m}{4BCD-1} \right)= (-1)^{(m-1)/2}(-1)^{(m-1)/2}=1\]
\end{dem}

\subsection{Modular equations }

For greater convenience, we call \textit{modular equation} a modular equation (or a system of modular equations) \emph{with constant coefficients}.\\

Since $A$ and $B$ play symmetrical roles, we may suppose%
\footnote{The arbitrary definition of $F$ ($A$ is factored out rather than $B$) was made in anticipation of this relation. Otherwise we could not have $\dg F=0$.}
that $\dg B\leqslant \dg A$, where $\dg$ is the degree of a polynomial.

\begin{lemme} \label{deg}
Let $p$ be a prime polynomial of degree 1.
\begin{enum}
\item[i)] If the relation (\ref{Rosati_1}) $ 4ABCD=A+B+pC$ holds, then 
\begin{subequations}\label{deg1}
\begin{align}
\dg A=1 \qquad \dg B=0 \qquad \dg C=0 \qquad \dg D=0 \qquad \label{deg1a}\\
\dg A=0 \qquad \dg B=0 \qquad \dg C=0 \qquad \dg D=1 \qquad \label{deg1b}\\
\dg A=1 \qquad \dg B=0 \qquad \dg C=1 \qquad \dg D=0 \qquad \label{deg1c}
\end{align} 
\end{subequations} 

\item[ii)] If the relation (\ref{Rosati_2}) $4ABCD=p(A+B)+C$ holds, then 
\begin{subequations}\label{deg2}
\begin{align}
\dg A=0 \qquad \dg B=0 \qquad \dg C=0 \qquad \dg D=1  \qquad \label{deg2a}\\
\dg A=1 \qquad \dg B=0 \qquad \dg C=0 \qquad \dg D=1  \qquad \label{deg2b}\\
\dg A=1 \qquad \dg B=0 \qquad \dg C=1 \qquad \dg D=0  \qquad \label{deg2c}\\
\dg A=2 \qquad \dg B=1 \qquad \dg C=0 \qquad \dg D=0  \qquad \label{deg2d}
\end{align} 
\end{subequations}
\end{enum}   
\end{lemme}
\begin{dem}
Since $\dg B\leqslant \dg A$ then $\dg (A+B)=\dg A$. Hence, by $C\,E=A+B$, we have $\dg C\leqslant \dg A$\\
 
i) By (\ref{Rosati_1}) it follows 
\begin{equation}
     (4ABD-p)C=A+B  \label{P1}  
\end{equation}
and
\begin{equation}
       (4BCD-1)A=B+Cp  \label{P2} 
\end{equation}

By (\ref{P1}) we have $\dg (4ABD-p)\leqslant \dg A$ and hence  $\dg (4ABD-p)\leqslant \dg ABD$.\\ 

* Case $\dg (4ABD-p)=\dg ABD$ \\
By (\ref{P1}) we have  \[ \dg ABD+\dg C=\dg A \]
and then 
    \[ \dg B+\dg C+\dg D=0  \label{P1.1} \]
This result implies, in view of (\ref{P2}), that 
            \[\dg A=\dg p=1\]

* Case $\dg (4ABD-p)<\dg ABD$. In this case\\
  \[ \dg ABD=\dg p=1  \et   \dg (4ABD-p)=0 \]
We have, by the first equation
           \[ \dg B=0 \et \dg A+\dg D=1 \] 
and by the second  together with (\ref{P1}) 
  \[ \dg C=\dg A \]

ii) By (\ref{Rosati_2}) it follows 
\begin{equation}
     (4ABD-1)C=p(A+B)  \label{P3}  
\end{equation}
and
\begin{equation}
       (4BCD-p)A=pB+C  \label{P4} 
\end{equation}

By (\ref{P3}) we have   
       \[ \dg ABD+\dg C=\dg A+\dg p \]
and then
       \[ \dg B+\dg C+\dg D=1 \]
Here $F=4BCD-p$. Then
   \[ \dg F \leqslant1 \et \dg (pB+C)=\dg pB\]
and together with (\ref{P4})
   \[  \dg F+\dg A=\dg B+1 \]     

* Case $\dg F= 1$. In this case $\dg A=\dg B$. On the other hand, as $\dg F=1$, there exists $t_0\in\R$ such that $F(t_0)=0$. From $FE=4B^2D+1$, it follows $4B^2(t_0)D(t_0)+1=0$ and then $D(t_0)<0$. Therefore $\dg D=1$ and $\dg B=\dg C=0$.\\

* Case $\dg F= 0$. In this case $\dg A=\dg B+1 $.    
\end{dem}

\begin{prop} \label{eqmod}
Let $p$ be a prime polynomial of degree 1. The fraction $4/p$ is 3-Egyptian if and only if one of the next 7 modular equations  holds.
\begin{subequations}\label{eqmod1}
\begin{align}
\hspace{1cm}& B+pC=0 \mod 4BCD-1   \label{eqmod1a} \\
\hspace{1cm}& p+E=0 \mod 4AB \et A+B=0 \mod E \qquad\label{eqmod1b} \\
\hspace{1cm}& p+E+4B^2D=0 \mod 4BDE   \label{eqmod1c}  
\end{align}
\end{subequations}
\begin{subequations}\label{eqmod2}
\begin{align}
\hspace{1cm}&  pE+1=0 \mod 4AB \et A+B=0 \mod E  \label{eqmod2a} \\
\hspace{1cm}&  p+F=0 \mod 4BC \et pB+C=0 \mod F   \label{eqmod2b}  \\
\hspace{1cm}&  p+F=0 \mod 4BD \et 4B^2D+1=0 \mod F  \label{eqmod2c} \\
\hspace{1cm}&   p+F=0 \mod 4CD \et p^2+4C^2D=0 \mod F \label{eqmod2d}
\end{align}
\end{subequations}
where $(A,B)=(B,C)=(C,D)=(4ABD,E)=(4BCD,F)=1$.
\end{prop}

\begin{dem} $\quad$ The $[\:]$ refer to the equations of the Lemma.
\begin{itemize}
\item[(\ref{eqmod1})] Here "(\ref{Rosati_eq}) is equivalent to (\ref{Rosati_E})" is written
	 \[ 4ABCD=A+B+pC \iff (p+E=4ABD \et A+B=CE)\]

\begin{itemize}[leftmargin=1ex]
\item[{(\ref{eqmod1a})}] Case [\ref{deg1a}] : $B,C,D$ are constants. If we suppose that (\ref{Rosati_eq}) holds, then
  \[  B+pC=(4BCD-1)A=0 \mod (4BCD-1) \] 
Conversely, we set 
  \[ A=\dfrac{B+pC}{4BCD-1} \]

\item[(\ref{eqmod1b})] Case [\ref{deg1b}] : $A,B,E$ are constants. If we suppose that (\ref{Rosati_eq}) holds, then 
  \[  p+E=4ABD=0 \mod 4AB \et A+B=CE=0 \mod E \] 
Conversely, we set 
  \[ D=\dfrac{p+E}{4AB} \et C=\dfrac{A+B}{E} \] 

\item[(\ref{eqmod1c})] Case [\ref{deg1c}] : $B,D,E$ are constants.
If we suppose that (\ref{Rosati_eq}) holds, then  
	\[ p+E=4(CE-B)BD \]
Hence   
  \[ p+E+4B^2D=4BDEC=0 \mod 4BDE \]  
Conversely, we set 
  \[  A=\dfrac{p+E}{4BD} \et C=\dfrac{p+E+4B^2D}{4BDE} \quad
  \left( C E = A+B\right) \]           
\end{itemize}

\item[(\ref{eqmod2})] Here "(\ref{Rosati_eq}) is equivalent to (\ref{Rosati_E})" is written
 \[ 4ABCD=p(A+B)+C \iff (4ABD=pE+1 \et A+B=CE) \]
and "(\ref{Rosati_eq}) is equivalent to (\ref{Rosati_F})" is written
\[ 4ABCD=p(A+B)+C \iff (p+F=4BCD \et pB+C=FA) \]
where $ F=4BCD-p \et FE=4B^2D+1 $
\begin{itemize}[leftmargin=1ex]
\item[(\ref{eqmod2a})] Case [\ref{deg2a}] : $A,B,E$ are constants. If we suppose that (\ref{Rosati_eq}) holds, then 
  \[ pE+1=4ABD=0 \mod 4AB \et A+B=CE=0 \mod E \]
Conversely, we set 
  \[ D=\dfrac{pE+1}{4AB} \et C=\dfrac{A+B}{E} \]

\item[\textit{In}]\textit{the next cases $\dg A=\dg B+1$ and then $\dg F=0$}.\\
 
\item[(\ref{eqmod2b})] Case [\ref{deg2b}] : $B,C,F$ are constants. If we suppose that (\ref{Rosati_eq}) holds, then  
  \[ p+F=4BCD=0 \mod 4BC \et pB+C=FA=0 \mod F \]
Conversely, we set 
  \[ D=\dfrac{p+F}{4BC} \et A=\dfrac{pB+C}{F} \]       
%
\item[(\ref{eqmod2c})] Case [\ref{deg2c}] : $B,D,F$ are constants. If we suppose that (\ref{Rosati_eq}) holds, then 
  \[ p+F=4BCD=0 \mod 4BD \et 4B^2D+1=EF=0 \mod F \]
Conversely, we set 
  \[ C=\dfrac{p+F}{4BD} \quad E=\dfrac{4B^2D+1}{F} \et A=CE-B\]
We observe that $FA=pB+C$.    
%
\item[(\ref{eqmod2d})] Case [\ref{deg2d}] : $C,D,F$ are constants. If we suppose that (\ref{Rosati_eq}) holds, then 
    \[ p+F=4BCD=0 \mod 4CD \et pB+C=FA=0 \mod F \]
As
	\[ pB+C=p\dfrac{p+F}{4CD}+C=\dfrac{p^2+pF+4C^2D}{4CD} \]    
it follows, since $(4CD,F)=1$, 
  \[ p^2+4C^2D=0 \mod F \] 
Conversely, we set 
  \[ B=\dfrac{p+F}{4CD} \et A=\dfrac{pB+C}{F} \]               
\end{itemize}

\end{itemize}
\end{dem}

We observe that, if $p$ is a prime polynomial of degree 1, the Lemme\ref{deg} shows that there are only 7 distinct cases, according to the degree of $A$, $B$, $C$, $D$ ($\dg B \leqslant\dg A $). By the Proposition~\ref{eqmod}, each case is connected to a modular equation. Hence, there exist only 7 distinct modular equations with \textit{constant coefficients}. So, we can build an algorithm giving the set (maybe empty) of \textit{all} the way to write $4/p$.

\subsection{Application to the integers}

The proof of the Proposition~\ref{eqmod} $\,$gives us formulas for $A,B,C,D$. These variables take strictly positive values when the given data are strictly positive and one of the equation (\ref{Rosati_1}) or (\ref{Rosati_2}) holds. Hence we have the following corollary.

\begin{corolaire}
Let $p$ be an odd prime integer. The fraction $4/p$ is 3-Egyptian if and only if one of the 7 modular equations of the Proposition~\ref{eqmod} holds.
\end{corolaire} 

Thereafter, we call these equations \textit{reference equations} not only for the polynomials but for the integers too. 

\subsubsection*{Comparison with previous results}

Four of these equations have been well known for a long time, but the others are new.

\begin{enum}
\item  Rosati \cite{Rosati} (1954) gives only one condition for (\ref{Rosati_1}) and one for (\ref{Rosati_2}). Although they are not written in a modular form, his equations (3) and (6)  are equivalent to (\ref{eqmod1a}) and (\ref{eqmod2a}).

\item Yamamoto \cite{Yamamoto} (1965) gives two conditions for (\ref{Rosati_1}) and two for (\ref{Rosati_2}). Written in a modular  form, his equations (3) to (6) are equivalent (not in the same order) to (\ref{eqmod1a}), (\ref{eqmod1b}), (\ref{eqmod2a}), (\ref{eqmod2b}).   
\end{enum}

Polynomials explain why the Yamamoto equivalent equations give distinct results. Even better, they give us three new equations. 

\subsubsection*{"Complete" set of modular equations}

Regarding prime polynomials of degree 1, the 7 reference equations  form a complete set%
\footnote{Moreover, example~2 below shows that these equations are independent.},
that is, if a modular equation $n=b \mod a$ (where $(a,b)=1$) is not equivalent to one of the reference equations then $4/(at+b)$ cannot be an 3-Egyptian fraction. This feature does not hold for integers : it may exist a process using such an equation and leading to the conclusion that $4/n$ is a 3-Egyptian fraction. But, in this case, this process  has to be of a still unknown new type.

\subsection{Examples}

\textbf{Example 0.} Of course, we may find the identities of paragraph \ref{reduction}. Here, we don't look after \emph{all} the way to write $4/p$, just those given in the paragraph.

\begin{itemize}
\item $p=3t-1$ verifies (\ref{eqmod1a}) : $p+1=0 \mod 3$ \\
where $B=C=D=1$, and hence $A = (p+1)/3=t$.
\item $p=4t-1$ verifies (\ref{eqmod1b}) : $p+1 = 0 \mod 4$ \\
where $A=B=E=1$ and hence $C=2$ et $D=(p+1)/4=t$. 
\item $p=8t-3$ verifies (\ref{eqmod1b}) : $p+3=0 \mod 8$ \\
where $A=1$, $B=2$, $E=3$ and hence $C=1$, $D=(p+3)/8=t$.
\end{itemize} 
  
\textbf{Example 1.} $\quad p=24\cdot 5t-23 \qquad (p=1 \mod 24 \et p=2 \mod 5)$\\
We give all the way to write $4/p\,$ and the distinctive feature is that the 7 reference equations (shown in $[\;]$) are present. We don't know another analogous example where $p=1 \mod 24$.
\begin{enum}
\item[{[\ref{eqmod1a}]}] $\quad \dfrac{4}{p}=\dfrac{1}{p}\left(\dfrac{1}{4}+\dfrac{1}{4(16t-3)} \right)+\dfrac{1}{2(16t-3)}$\\

\item[{[\ref{eqmod1b}]}] $\quad \dfrac{4}{p}=\dfrac{1}{p}\left(\dfrac{1}{10(6t-1)}+\dfrac{1}{2(6t-1)} \right)+\dfrac{1}{5(6t-1)}$\\

\item[{[\ref{eqmod1c}]}] $\quad \dfrac{4}{p}=\dfrac{1}{p}\left(\dfrac{1}{10t}+\dfrac{1}{10t(6t-1)} \right)+\dfrac{1}{5(6t-1)}$\\

                         $\quad \dfrac{4}{p}=\dfrac{1}{p}\left(\dfrac{1}{2t}+\dfrac{1}{2t(15t-1)} \right)+\dfrac{1}{2(15t-1)}$\\

\item[{[\ref{eqmod2a}]}] $\quad \dfrac{4}{p}=\dfrac{1}{5(21t-4)}+\dfrac{1}{2(21t-4)} +\dfrac{1}{10(21t-4)\,p}$\\

\item[{[\ref{eqmod2b}]}] $\quad \dfrac{4}{p}=\dfrac{1}{5(6t-1)}+\dfrac{1}{2(6t-1)(100t-19)} +\dfrac{1}{10(6t-1)(100t-19)\,p}$\\

                     $\quad \dfrac{4}{p}=\dfrac{1}{5(6t-1)}+\dfrac{1}{10(6t-1)(20t-3)} +\dfrac{1}{2(6t-1)(20t-3)\,p}$\\

                     $\quad \dfrac{4}{p}=\dfrac{1}{2(15t-1))}+\dfrac{1}{(15t-1)(16t-3)} +\dfrac{1}{2(15t-1)(16t-3)\,p}$\\
              
\item[{[\ref{eqmod2c}]}] $\quad \dfrac{4}{p}=\dfrac{1}{5(6t-1)}+\dfrac{1}{10(6t-1)(21t-4)} +\dfrac{1}{10(21t-4)\,p}$\\

\item[{[\ref{eqmod2d}]}] $\quad \dfrac{4}{p}=\dfrac{1}{5(6t-1)}+\dfrac{1}{10(120t^2-43t+4)} +\dfrac{1}{10(6t-1)(120t^2-43t+4)\,p}$\\
\end{enum} 

\textbf{Example 2.}$\quad$ In this example, each $p\,$ is of the form $p=24\cdot 583t+b$.  At the opposite of the example 1, the distinctive feature is that, for some $b$, there is only one way to write $4/p$. A value of $b$ is given for each reference equation.
\begin{enum}

\item[{[\ref{eqmod1a}]}] $\quad  p=24\cdot 583t-911 \quad  (p=1 \mod 24 \et p=255 \mod 583)$\\

      $\quad \dfrac{4}{p}=\dfrac{1}{p}\left( \dfrac{1}{2\cdot 73} +\dfrac{1}{6(16t-1)} \right)+\dfrac{1}{3\cdot 73(16t-1)}$\\

\item[{[\ref{eqmod1b}]}] $\quad p=24\cdot 583t-119 \quad  (p=1 \mod 24 \et p=464 \mod 583)$\\ 

      $\quad \dfrac{4}{p}=\dfrac{1}{p}\left( \dfrac{1}{66t} +\dfrac{1}{53t} \right)+\dfrac{1}{66\cdot 53t} $\\
  
\item[{[\ref{eqmod1c}]}] $\quad p=24\cdot 583t-1127 \quad  (p=1 \mod 24 \et p=39 \mod 583)$\\

     $\quad \dfrac{4}{p}= \dfrac{1}{p}\left(\dfrac{1}{22t} +\dfrac{1}{ 2t(159t-11)}\right) +\dfrac{1}{ 22(159t-11)} $\\
  
\item[{[\ref{eqmod2a}]}] $\quad p=24\cdot 583t-1799 \quad  (p=1 \mod 24 \et p=533 \mod 583)$\\ 

     $\quad \dfrac{4}{p}= \dfrac{1}{50\cdot 1749(70t-9)} +\dfrac{1}{50(70t-9)} +\dfrac{1}{1749(70t-9)\,p}$\\
 
\item[{[\ref{eqmod2b}]}] $\quad p=24\cdot 583t-11159 \quad  (p=1 \mod 24 \et p=501 \mod 583)$
\begin{multline*}
     \dfrac{4}{p}= \dfrac{1}{ 22(159t-125)}
     +\dfrac{1}{ 8(242t-193)(159t-125)}\\
     +\dfrac{1}{88(242t-193)(159t-125)\,p}
\end{multline*}
  
\item[{[\ref{eqmod2c}]}] $\quad p=24\cdot 583t-503 \quad  (p=1 \mod 24 \et p=80 \mod 583)$\\

     $\quad \dfrac{4}{p}= \dfrac{1}{6\cdot 583t} +\dfrac{1}{ 6\cdot 53t(306t-11)} +\dfrac{1}{583(306t-11)\,p} $\\
 
\item[{[\ref{eqmod2d}]}] $\quad p=24\cdot 583t-6407 \quad  (p=1 \mod 24 \et p=6 \mod 583)$
\begin{multline*}
     \dfrac{4}{p}= \dfrac{1}{22(159t-71)} 
     +\dfrac{1}{22(13992t^2-12655t+2861)}\\ 
       +\dfrac{1}{11(159t-71)(13992t^2-12655t+2861)\,p} 
\end{multline*} 

\end{enum}

\section{Modular sieve}

The algorithms setting, for a given integer $n>2$, at least one way (and even more) to write $4/n$ are interesting. However, regarding the checking of the conjecture, an efficient algorithm needs an another point of view%
\footnote{This point of view has been used since Rosati's paper (or maybe before).}.\\

We denote by $\N_0$ the set of the integers $n\in\N$ verifying the condition $n=1 \mod 24$. The process described below takes account specifically of the fact that the checked integers are in $\N_0$. On the other hand, we let down the condition that $n$ is prime, which needs too much running time. Regarding the polynomial $at+b$, the correlated conditions are $at+b =1 \mod 24$ (which is equivalent to $a=0 \mod 24$ and $b=1 \mod24$) and the cancellation of the condition $(a,b)=1$.

\subsection{Modular filters}
\textit{Definition : } A \textit{sieve} is a sorted set of filters.\\

\textit{Definition : } A \textit{filter}%
\footnote{We use the terminology given by Swett. If an integer $n\in\N_0$ is such that $n\% m\in F$ then $n$ verifies the conjecture and $n$ is "trapped" by the filter. Otherwise $n$ "pass through".}%
\textit{modulo} $m$ is a set $F$ such that for any $n\in\N_0$    
\[ n\%m\in F \Rightarrow 4/n \;\textnormal{is 3-Egyptian}\]
where  $n\%m$ is the residue of $n$ modulo $m$ (notation borrowed from C language).\\

For $a>0$, we denote by $\Omega_a$ the set of $b\in\Z$ such that 
$4/(at+b)$ is a 3-Egyptian fraction. If $m$ is odd, we set 
$ S_m=\left(\Omega_{[m,24]}\cap\N_0\right)\%m $ 
where $[u,v]=\ppcm(u,v)$. It follows some obvious proprieties.

\begin{enumerate}[label =\textit{\roman*)}]
\item  if $q\mid a$ then $\Omega_q\subset\Omega_a$.
\item if $b_1=b_2 \mod a$ then 
$b_1\in\Omega_a \Rightarrow b_2\in\Omega_a$.
\item if $n\in \Omega_a$ ($n>0$) then $4/n$ is a 3-Egyptian fraction.
\item if $n\in\N_0$ then $\left( n\in\Omega_{[m,24]} \iff n\%m\in S_m \right)$, which shows that $S_m$ is a filter modulo $m$.
\item \label{reduit} if $n\in\N_0$ and if $q \mid m$ then $n\%q \in S_q \Rightarrow n\%m \in S_m$. 
\end{enumerate}

\textit{Définition :} We say that $n\in\N_0$ is \emph{certified} if there exists $m$ such that $n\%m\in S_m$. We also say that $n$ is certified by $m$ or that $m$ is a \emph{modular certificate} of $n$ (vocabulary borrowed from the complexity theory).\\ 

\textit{The first results with prime integers} :
\begin{multicols}{2}
\begin{enum}
\item[] $S_5=\{0,2,3\}$
\item[] $S_7=\{0,3,5,6\}$
\item[] $S_{11}=\{0,7,8,10\}$
\item[] $S_{13}=\{0,5,6,8,11\}$
\item[] $S_{17}=\{0,10,11,12,14\}$
\item[] $S_{19}=\{0,8,12,14,15,18\}$
\item[] $S_{23}=\{0,7,10,11,15,17,19,20,21,22\}$
\item[] $S_{29}=\{0,14,18,19,21,26,27\}$
\item[] $S_{31}=\{0,15,22,23,24,27,29,30\}$
\item[] $S_{37}=\{0,5,15,18,22,23,29,32,35\}$
\end{enum}
\end{multicols}

\textit{Some results with odd composite integers} : 
\begin{enum}
\item[]$S_{15}=\{7, 10, 13\}$
\item[]$S_{35}=\{0,2,3,5,6,7,8,10,12,13,14,15,17,18,19,20,21,22,23,24 $\\ \hphantom{m} \hfill $,25, 26,27,28,30,31,32,33,34\}$
\item[]$S_{55}=\{0,2,3,5,7,8,10,11,12,13,15,17,18,19,20,21,22,23 $\\ \hphantom{m}\hfil $,24,25,27,28,29,30,32,33,35,37,38,39,40,41,42,43$\\ \hphantom{m}\hfill $,44,45,47,48,50,51,52,53,54\} $

\end{enum}

\subsection{Shortened filters}   
If $m$ is composite, some integers $n\in\N_0$ are certified both par $m$ and by one of its divisors (cf. the propriety~\ref{reduit} above). The next definition allows us to point out what is particular to $m$.\\

\textit{Definition} : The \textit{shortened filter} $S^*_m$ is the set of all $x\in S_m$ such that
\begin{center}
	$x\%q\notin S_q$ for any $q\mid m$, $q\neq m$
\end{center}

We observe that if $m$ is prime then $S^*_m=S_m$.\\

\textit{The first (no empty) results}

\begin{multicols}{2}
\begin{enum}
\item[] $S^*_{55}=\{24,39\}$ 
\item[] $S^*_{65}=\{54,59\}$ 
\item[] $S^*_{77}=\{46,72\}$
\item[] $S^*_{85}=\{54,74\}$ 
\item[] $S^*_{95}=\{29,59,79,89\}$ 
\item[] $S^*_{99}=\{61,79,94\}$ 
\item[] $S^*_{117}=\{85,106\}$
\item[] $S^*_{119}=\{23,39,57,58,71,88,107,109\}$
\end{enum}
\end{multicols}

\section{Checking of the conjecture}

\subsection{Choice of the progressions}

The checked integers $n$ are in an arithmetic progression, namely they are of the form $n=24k+1$. We call \textit{gap} of the progression the difference between two consecutive terms. Here the gap is $G_0=24$ but if we use some filters $S_m$ we may obtain other progressions whose gap is bigger.\\

With $S_5=\{0,2,3\}$ we check only $n$ such that 
\[ n\%24 = 1 \et n\%5\in \{1,4\} \] 
and hence, by the Chinese remainder theorem
\[ n\%120\in \{1,49\} \]
The new gap is $G_1=120$, and there are 2 residues : then the mean gap is $g_1=60$. In comparison to 24, we check $2.5$ times fewer integers ($60/24=2.5$).\\

Next, with  $S_7=\{0,3,5,6\}$ we check only $n$ such that
\[ n\%120 \in \{1,49\} \et n\%7 \in \{1,2,4\} \] 
and hence
\[ n\%840 \in R_2 \]
where $R_2=\{1,121,169,289,361,529\}$ is the set of the residues%
\footnote{It was the choice made by Swett.}.
The new gap is $G_2=840$ and the mean gap is $g_2=140$. In comparison to 24, we check nearly 6 times fewer integers ($140/24=35/6$).\\
 
We may keep on and use others $S_m$. The checked integers are then of the form  
\[ n\%G_i\in R_i \]
where the first values of $G_i=G_{i-1}m_i$ and $\#R_i$ (the number of elements of $R_i$) are set out in the following table.

\[\begin{array}{|c|c|c|c|c|}
\Hline i & m_i &       G_i     & \#R_i &   g_i \\ 
\Hline 1 & 5 &           120 &         2 &   60 \\
\Hline 2 & 7 &           840 &         6 &   140 \\ 
\Hline 3 & 11 &       9\,240 &        34 &   272 \\ 
\Hline 4 & 13 &      120\,120 &       192 &   626 \\ 
\Hline 5 & 17 &   2\,042\,040 &    1\,507 &  1\,355 \\ 
\Hline 6 & 19 & 38\,798\,760 &   13\,380 &  2\,900 \\
\Hline 7 & 23 & 892\,371\,480 &   147\,348 &  6\,056 \\ 
\hline 
\end{array}\]\\ 

Three comments about this table. 
\begin{itemize}
\item The first concerns the reduction of $R_i$ (done in the table). If $n=r \mod G_i$ then for any $q$ divisor of $G_i$ we have 
$n\%q = r\%q$. Hence, we may remove the residues $r\in R_i$ verifying $r\%q \in S_q$. This reduction is essential, otherwise it's just a useless complicated process. 

\item The second concerns the last column : the mean gap $g_i=G_i/\#R_i$ is a good speed indicator. By example, as $6\,056/140\approx 43$, then using $G_7$ rather than $G_2$ leads to check about $43$ times fewer integers and the running time is shortened accordingly. 

\item The last concerns the choice of the $m_i$. The usual order is misleading~: each other set of seven integers seems to give a worse $g_7$. Next, with height integers we expect to add $31$ (rather than $29$). However, these two propositions have to be confirmed.
\end{itemize}

\subsection{Optimized sieve}
We denote by $\N_i$ the set of all the integers $n\in\N$ verifying $n\%G_i\in R_i$. As the conjecture is verified for any integer $n\notin N_i$, we have just to check the prime integers of $N_i$.\\
  
Let $N=10^{17}$ and $M$ the set of all the odd integers $m<5\,000$. We claim that each $n\in\N_7$ has a modular certificate in $M$ if $n<N$ and if $n$ is not a square. It is equivalent to say that $ \N_7 \setminus \bigcup_{m\in M}\Omega_{[m,24]}$ has not any element $n<N$, except squares.\\

We could use this $M$ to prove that the conjecture is verified up to $N$. However, if we want a running time as fewer as possible, we have to optimize the sieve. For this purpose, we remove the useless elements and sort $M$ in order to have at first the most efficient filters%
\footnote{The approach mostly hinge on experiments and make use of the shortened filters.}.
By example for  $N=10^{17}$, we give below the set $M=MOD$ which is used in our C++ program.\\

$ MOD =\{$ 3, 5, 7, 11, 13, 17, 19, 23, 4495, 2491, 2627, 4661, 4223, 1505, 4355, 3355, 4509, 4775, 2629, 4565, 4599, 4585, 3955, 3535, 3857, 3115, 3419, 3949, 3395, 3353, 1391, 1199, 3775, 4325, 4031, 2799, 1639, 4475, 2159, 4795, 2961, 1727, 4075, 1791, 4743, 2849, 3595, 1115, 3445, 3263, 2155, 2065, 2515, 2681, 4195, 3223, 2519, 4103, 3731, 4345, 3743, 2439, 1055, 2951, 1799, 4193, 1991, 3047, 2933, 3951, 4147, 1631, 2219, 4615, 3913, 3679, 1535, 2959, 1655, 4123, 1439, 3839, 1319, 3695, 4255, 3895, 1351, 2495, 1835, 2855, 2335, 4529, 1917, 1079, 1559, 1735, 1679, 2165, 4367, 4555, 2359, 2723, 3065, 3899, 3295, 3035, 4927, 3359, 4437, 3635, 4315, 2735, 3241, 4319, 4105, 4069, 1039, 4059, 1247, 3095, 4571, 3665, 1007, 1583, 4895, 1847, 2435, 1765, 2807, 3647, 1343, 2651, 3965, 1511, 2655, 4403, 1151, 887, 2935, 3545, 2879, 1967, 2815, 2399, 4419, 1159, 4487, 3119, 1223, 2039, 4745, 2305, 1103, 4077, 3215, 3715, 2279, 4915, 4873, 1031, 1475, 3865, 2483, 1399, 1823, 3173, 3305, 2241, 3985, 3563, 1349, 1259, 3959, 4415, 3455, 2615, 1487, 3599, 3935, 1759, 3505, 1871, 4879, 4535, 3199, 2045, 1367, 1493, 1919, 3787, 2111, 1975, 2053, 4739, 1231, 4151, 1837, 1213, 3655, 2183, 4135, 4939, 1019, 3023, 3995, 1855, 4265, 4079, 3983, 2575, 1063, 2351, 4985, 2687, 3167, 2447, 2725, 4631, 4595, 4115, 4175, 4055, 4679, 1013, 2239, 4385, 1091, 3429, 1909, 1719, 2365, 3415, 3079, 4955, 1147, 1133, 3191, 3475, 2759, 4405, 2207, 4765, 3431, 1139, 4471, 2727, 4145, 3247, 1279, 1751, 3755, 1087, 4835, 1733, 4645, 1979, 4711, 1177, 1073, 3055, 3239, 2999, 2087, 4855, 4039, 1703, 3527, 4295, 4799, 4207, 4505, 1187, 1109, 1567, 1379, 2119, 2911, 2591, 2015, 3785, 1651, 3155, 1819, 4751, 3719, 4735, 2345, 2831, 2099, 4995, 1427, 2059, 1333, 1069, 1663, 2719, 2063, 4285, 2231, 1093, 1607, 1423, 1411, 1027, 3805, 1769, 1121, 1903, 4063, 4759, 1363, 1973, 4715, 2663, 3863, 1433, 2479, 4703, 3299, 1451, 2339, 1613, 1471, 1619, 3671, 2287, 2367, 3845, 3537, 1591, 3733, 4463, 1271, 1931, 4619, 2903, 2135, 4921, 4685, 4705, 1003, 1429, 1193, 4067, 3275, 4311, 1327, 3015, 1499, 2413, 1237, 1181, 4045, 4081, 3605, 3779, 3103, 2837, 1579, 3439, 1033, 3799, 2333, 1829, 1241, 4393, 2357, 4159, 2699, 3791, 2453, 3625, 2579, 4945, 4127, 1649, 4741, 4871, 1667, 2177, 3835, 1043, 3407, 4919, 4885, 2267, 2693, 2507, 4967, 2327, 4639, 1691, 1549, 2583, 1123, 1717, 1999, 1807, 1933, 4553, 1049, 3479, 1553, 1853, 2543, 4343, 1501, 2743, 3699, 1787, 3989, 1129, 1525, 4445, 1675, 1993, 1301, 2273, 1217, 1843, 4003, 2411, 3245, 3401, 1117, 1789, 3379, 3901, 1831, 1957, 4085, 1507, 1987, 3373, 3893, 1621, 1943, 3937, 1291, 1543, 1571, 2143, 2533, 2767, 3253, 4883, 2551, 2833, 1229, 1877, 1949, 2009, 4391, 1643, 2251, 2729, 3915, 1907, 2243, 2603, 2669, 2897, 3043, 3313, 3739, 1171, 1361, 1817, 1879, 2659, 3623, 4283, 4859, 1537, 2003, 2161, 2389, 2869, 4439, 1099, 1415, 2269, 2293, 2943, 3233, 3967, 4181, 4261, 4559, 4699, 1447, 1895, 1921, 2195, 2939, 3293, 3565, 3607, 3749, 4247, 4591, 4829, 1157, 1417, 1951, 1997, 2179, 2225, 2619, 2785, 3041, 3717, 4583, 4783, 4887, 1283, 1517, 1721, 1747, 1961, 2033, 2117, 2129, 2741, 2803, 2893, 3161, 3589, 3613, 1211, 1273, 1459, 1483, 1811, 1867, 1889, 1971, 2043, 2069, 2149, 2213, 2423, 2709, 2779, 3013, 3149, 3551, 4013, 4097, 4363, 4399, 4589, 4681, 1021, 1097, 1145, 1197, 1297, 1373, 1397, 1555, 1609, 1723, 1773, 1777, 1783, 1801, 2123, 2191, 2259, 2291, 2371, 2407, 2443, 2671, 2845, 3389, 3493, 3725, 4021, 4171, 4351, 4999 $\}$

\subsection{Results}
The checked integers are of the form $n=r+k\times G_7$ where $r\in R_7$ and $0 \leqslant k < K$. With $N=10^{17}$, we take $K=112\,066\,560$. Therefore we check $16\,512\,783\,482\,880$ integers including $51\,732\,427$ squares.\\

For each $m\in MOD$, the number of integers certified by $m$ is given at the same rank in the table below. We may observe that the sum of these numbers added with the number of squares is equal to the number of checked integers.\\

$ \{$ 0, 0, 0, 0, 0, 0, 0, 9223757362766, 3739609092281, 1565954748220, 739166512371, 397180210351, 249398230928, 169050837573, 104088377604, 69101085771, 54368854713, 42523071218, 33206924179, 23406992663, 18890746142, 15211918708, 11968966501, 9473482721, 7560449664, 6273004978, 5196086887, 4344239727, 3485872879, 2944121141, 2498890993, 2067185415, 1765012627, 1499112458, 1259328652, 1044404123, 874512654, 723079141, 617340453, 515245196, 452563855, 390773540, 343076561, 300065591, 260653549, 229207022, 198772233, 174906642, 153551008, 135203129, 118673167, 99017032, 88208940, 78571579, 69928806, 62430095, 55603526, 48999877, 42755472, 38618483, 34775913, 31335757, 27560576, 24702471, 22410294, 19685100, 17852098, 16081935, 14466854, 13141729, 11726132, 10640116, 9491430, 8477371, 7732328, 6982821, 6272905, 5702788, 5268793, 4722801, 4390120, 4019516, 3650026, 3398755, 3140726, 2945648, 2736821, 2552135, 2394011, 2241501, 2100950, 1967613, 1834764, 1721795, 1606462, 1497392, 1396075, 1313066, 1211933, 1135277, 1058550, 992002, 932632, 867721, 807226, 759519, 707804, 665956, 625412, 586324, 552400, 520156, 484882, 459951, 434799, 408981, 385709, 365271, 343865, 322175, 305617, 290089, 275265, 260444, 247858, 235278, 223480, 212702, 201616, 190609, 179783, 171406, 163172, 153416, 146567, 138229, 131604, 125502, 118787, 113170, 106706, 101331, 96415, 91243, 87184, 83366, 79429, 75744, 72027, 68424, 65511, 62568, 59796, 57362, 54579, 52058, 49809, 47397, 45319, 43371, 41545, 39868, 38029, 36324, 34724, 33179, 31866, 30561, 29163, 27958, 26789, 25721, 24727, 23684, 22631, 21616, 20727, 19908, 19064, 18235, 17562, 16604, 15933, 15087, 14379, 13832, 13303, 12763, 12305, 11871, 11311, 10807, 10411, 9972, 9650, 9215, 8834, 8518, 8200, 7898, 7612, 7318, 7021, 6760, 6493, 6232, 5968, 5727, 5528, 5308, 5119, 4884, 4693, 4538, 4369, 4224, 4063, 3940, 3786, 3653, 3530, 3414, 3246, 3117, 3013, 2912, 2787, 2685, 2559, 2455, 2384, 2269, 2185, 2103, 2048, 1985, 1908, 1831, 1777, 1725, 1671, 1621, 1563, 1517, 1465, 1423, 1375, 1335, 1285, 1245, 1196, 1156, 1122, 1076, 1048, 1012, 977, 943, 920, 892, 860, 834, 803, 783, 761, 743, 718, 701, 679, 654, 637, 602, 586, 571, 559, 539, 521, 510, 488, 473, 461, 445, 435, 426, 414, 402, 383, 374, 362, 352, 346, 335, 324, 317, 313, 303, 295, 284, 277, 270, 261, 252, 245, 242, 232, 226, 223, 217, 213, 205, 200, 197, 191, 187, 180, 177, 174, 169, 166, 162, 157, 153, 151, 147, 142, 138, 135, 133, 131, 129, 127, 124, 119, 115, 114, 111, 110, 105, 102, 100, 98, 95, 94, 92, 90, 86, 85, 83, 82, 79, 78, 78, 76, 75, 72, 70, 69, 67, 65, 64, 61, 60, 58, 57, 57, 55, 55, 54, 52, 51, 49, 48, 47, 46, 44, 44, 43, 42, 42, 41, 40, 39, 39, 38, 37, 37, 35, 35, 34, 34, 33, 32, 32, 30, 30, 30, 29, 29, 29, 27, 27, 26, 25, 25, 25, 24, 24, 24, 23, 23, 22, 22, 22, 21, 21, 20, 20, 19, 19, 19, 18, 18, 18, 17, 17, 17, 17, 16, 16, 16, 15, 15, 15, 15, 14, 14, 14, 13, 13, 13, 13, 13, 13, 13, 13, 12, 12, 11, 11, 11, 11, 11, 10, 10, 10, 10, 9, 9, 9, 9, 9, 9, 9, 9, 8, 8, 8, 8, 8, 8, 8, 8, 7, 7, 7, 7, 7, 7, 6, 6, 6, 6, 6, 6, 6, 6, 6, 6, 6, 5, 5, 5, 5, 5, 5, 5, 5, 5, 5, 5, 5, 4, 4, 4, 4, 4, 4, 4, 4, 4, 4, 4, 4, 4, 3, 3, 3, 3, 3, 3, 3, 3, 3, 3, 3, 3, 3, 3, 2, 2, 2, 2, 2, 2, 2, 2, 2, 2, 2, 2, 2, 2, 2, 2, 2, 2, 2, 2, 2, 2, 2, 2, 1, 1, 1, 1, 1, 1, 1, 1, 1, 1, 1, 1, 1, 1, 1, 1, 1, 1, 1, 1, 1, 1, 1, 1, 1, 1, 1, 1, 1, 1 $\}$



\begin{thebibliography}{99}
\bibitem {Bernstein} \textsc{Bernstein} Von Leon (1962), \textit{Zur Lösung der diophantischen Gleichung $m/n=1/x+1/y+1/z$, insbesondere im Fall $m=4$}, Journal für die reine und angewandte Mathematik, volume 211, p. 1-10.
\verb+http://gdz.sub.uni-goettingen.de/dms/load/img/?PPN=GDZPPN002179792+

\bibitem {Mordell} \textsc{Mordell} Louis Joel(1967), \textit{Diophantine Equations}, Academic Press, p. 287-290.   

\bibitem {Rosati} \textsc{Rosati} Luigi Antonio (1954), \textit{Sull'equazione diofantea } $4/n=1/x_1+1/x_2+1/x_3$, Bollettino dell'Unone Matematica Italiana, serie 3, volume 9 n.1 p. 59-63.
\verb+http://www.bdim.eu/item?fmt=pdf&id=BUMI_1954_3_9_1_59_0+  

\bibitem {Schinzel} \textsc{Schinzel} Andrzej (2000), \textit{On sums of three unit fractions with polynomial denominators}, Functiones et Approximatio Commentarii Mathematici vol.XXVIII p.187-194. \verb+http://www.staff.amu.edu.pl/~fa/XXVIII/fa-28-1-187.pdf+   

\bibitem{Swett} \textsc{Swett} Allan (1999), \textit{The \Erdos\ conjecture}, Current Research on ESC, rev.10/28/99. \verb+http://math.uindy.edu/swett/esc.htm+

\bibitem{Yamamoto} \textsc{Yamamoto} Koichi (1965), \textit{On the diophantine equation}  $4/n=1/x+1/y+1/z$, Memoirs of the Faculty of Science, Kyushu University, Series A, Mathematics, Vol. 19, p. 37-47. 
\verb+https://www.jstage.jst.go.jp/article/kyushumfs/19/1/19_1_37/_pdf+


\end{thebibliography}
\end{document}